%% file: lag5.tex
\documentstyle[psfig]{article}

\input mssymb12.tex
\input newmacro.tex

\title{Lagrangian spheres in $S^2\times S^2$}
\author{R. Hind\thanks{Supported in part by NSF grant DMS-0204634.}}

\date{\today}

\begin{document}

\maketitle

\section{Introduction}

The purpose of this paper is to consider the symplectic manifold
$S^2 \times S^2$ with the direct sum symplectic form
$\omega=\omega_0 \oplus \omega_0$, where $\omega_0$ is the
standard area form on $S^2$. Then the antidiagonal
$\overline{\Delta}=\{(x,-x)|x\in S^2\}\subset S^2\times S^2$ is a
Lagrangian submanifold, that is, $\omega|_L=0$. It is the aim of
this paper to demonstrate that any Lagrangian sphere $L$ in $S^2
\times S^2$ must actually be isotopic to $\overline{\Delta}$
through Lagrangian spheres (or, in other words, Lagrangian
isotopic to $\overline{\Delta}$). That is, we will show the
following.

\begin{theorem} Let $L$ be a Lagrangian sphere in $S^2 \times S^2$. Then there exists a Hamiltonian diffeomorphism of $S^2 \times S^2$ mapping
$L$ onto $\overline{\Delta}$.
\end{theorem}

We note that any such Lagrangian sphere in $S^2 \times S^2$ is
certainly homologous to $\overline{\Delta}$ since, by A.
Weinstein's Lagrangian neighborhood theorem, all Lagrangian
spheres have self-intersection number $-2$ .

The problem of studying Lagrangian knots in symplectic manifolds
was first proposed by V. I. Arnold in \cite{arnold} and some
results are already known in this direction, see for example the
survey \cite{survey}. We observe that $\Delta=\{(x,x)|x\in
S^2\}\subset S^2\times S^2$ is a symplectic submanifold and
$S^2\times S^2\setminus \Delta$ is symplectomorphic to a
neighbourhood of the zero-section in $T^* S^2$, mapping
$\overline{\Delta}$ to the zero-section. Y. Eliashberg and L.
Polterovich have shown in \cite{ep} that any Lagrangian
submanifold in $T^* S^2$ must be smoothly isotopic to the
zero-section. On the other hand, P. Seidel in \cite{seidel} has
given examples of Lagrangain spheres in symplectic manifolds which
are smoothly isotopic but not Lagrangian isotopic. A. Ivrii has
obtained some similar results to our own for Lagrangian tori.

The main point of our proof will be to construct two transverse
foliations of $S^2 \times S^2$ by spheres in the classes
$[\mathrm{point} \times S^2]$ and $[S^2 \times \mathrm{point}]$.
These spheres should be holomorphic with respect to some
almost-complex structure compatible with $\omega$ and each sphere
should intersect $L$ transversally in a single point. Of course,
without the condition of intersecting $L$ transversally we could
take the standard foliation by spheres which are holomorphic with
respect to a split complex structure $\Bbb C P^1 \times \Bbb C
P^1$. To obtain our required foliations we will start with this
complex structure but then deform it in a neighborhood of a
contact hypersurface. As the almost-complex structure is deformed,
Gromov showed in \cite{gr} that transverse foliations of
holomorphic curves will continue to exist, but they will also be
deformed and we will show that they eventually become transverse
to $L$. The deformation of the almost-complex structure that we
use is known as stretching-the-neck.

Taking a limit as the neck is stretched to infinite length is
possible by some recent results due to F. Bourgeois, Y.
Eliashberg, H. Hofer, K. Wysocki and E. Zehnder, see \cite{behwz},
one obtains finite energy holomorphic curves in symplectic
manifolds with cylindrical ends. One of these manifolds is $T^*
L=T^* S^2$ and we describe in section $3$ the resulting foliation
of this manifold by holomorphic curves. In section $4$ we discuss
the behaviour of holomorphic curves as we deform the
almost-complex structure (stretching the neck) and in section $5$
use the accumulated infomation to reach our conclusion on the
Lagrangian isotopy class. Various facts about finite energy
holomorphic curves in symplectic manifolds with cylindrical ends
and the relevant compactness theorem are gathered together first
in section $2$.

It is the author's pleasure to thank Yasha Eliashberg for various
enlightening conversations on these topics and also to thank the
I.A.S. for their hospitality while some of this work was done.

\section{Holomorphic curves in symplectic manfolds with cylindrical ends}

In this section we will state some theorems about holomorphic
curves in open symplectic manifolds with cylindrical ends. The
definitions and most of the proofs can be found in the series of
papers by H. Hofer, K. Wysocki and E. Zehnder, \cite{hofa},
\cite{hofi}, \cite{hoff}, \cite{hofd}. The generalizations to the
slightly degenerate situation which we will study are taken from
the paper of F. Bourgeois, \cite{bor}, see also \cite{bort} for
the proofs. Such a theory of holomorphic curves forms the basis of
symplectic field theory, see \cite{egh}.

We are interested in symplectic manifolds $(W,\omega)$ with
noncompact ends symplectomorphic to either $((0,\infty)\times
M,d(e^t\alpha))$ or $((-\infty,0)\times M,d(e^t\alpha))$. Here $M$
is a contact $3$-manifold with $\alpha$ a contact form. Let $X$ be
the corresponding Reeb vectorfield (which is uniquely defined by
$X \rfloor d\alpha =0$ and $\alpha (X)=1$) and $\xi=\{\alpha
=0\}$. The two types of ends are called convex or concave
respectively.

Given such an $(W,\omega)$, we equip it with a compatible almost-complex
structure $J$. The compatibility condition means that $\omega(V,JV)>0$
for all non-zero tangent vectors $V$, and
at the point $(a,m)\in (0,\infty)\times M$ or $(-\infty,0)\times M$, the almost-complex
structure is defined by
\begin{equation}
J(a,m)(h,k)=(-\alpha(m)(k), J'(m)\pi k+hX(m))
\end{equation}
for $(h,k)\in T_{(a,m)}((-\infty,\infty)\times M)$, where $\pi :TM\to \xi$
denotes the projection along $X$ and $J'$ is a fixed complex structure on $\xi$.

Finite energy holomorphic spheres are defined as follows.

Suppose $u:S^2\setminus \Gamma \to W$ is a proper map, where $S^2\setminus \Gamma$ denotes the Riemann sphere minus a finite set $\Gamma$ of punctures.

The energy of $u$ can be
defined by
$$E(u)=sup_{\phi}\int_{\C}u^{*}\omega_{\phi}$$
where the supremum is taken over all smooth, increasing
functions $\phi :(-\infty,\infty)\to (0,2)$ such that
$\phi =e^t$ for $t$ close to $0$ and
$\omega_{\phi}$ is defined
to be $d(\phi \alpha)$ on $(0,\infty)\times M$ and $(-\infty,0)\times M$ and equal to
the original $\omega$ elsewhere.

A finite energy holomorphic curve is then defined to be a $J$-holomorphic map $u$ with $E(u)<\infty$.

Let $\eta_t$ be the flow of the Reeb vectorfield $X$ on $M$
associated to $\alpha$. Suppose that $x$ is a periodic orbit of
$X$ of period $T$. Since $\eta_t^* \alpha = \alpha$, the
differential $D\eta_T$ induces a linear map $L:\xi_{x(0)}\to
\xi_{x(0)}$, and $x$ is called nondegenerate if $L$ does not
contain $1$ in its spectrum. The contact forms we will use do not
have nondegenerate periodic orbits. In fact our contact manifold
$M$ will actually be foliated by periodic orbits, all with the
same period, of the corresponding Reeb vectorfield. This is a
special case of a Reeb flow of Morse-Bott type.

The following theorem is proven by Hofer, Wysocki and Zehnder
in \cite{hofa}, Theorems $1.2$ and $1.4$, in the nondegenerate case and
in \cite{hofd} in the Morse-Bott situation.

\begin{theorem}
Let $u$ be a finite energy holomorphic sphere from a Riemann surface with a noncompact end identified with $\Bbb C\setminus D$. We may assume that
$u(\Bbb C\setminus D)\subset (-\infty,0)\times M$ or $(0,\infty)\times M$.
Let $\overline{u}$ denote the projection of $u$ to $M$.
Then there exists a periodic orbit $x$ of the
Reeb vectorfield $X$, say of period $T$, and a
sequence $R_k \to \infty$ such that
$$\overline{u}(R_k e^{\pm 2\pi it/T})\to x(t)$$ in $C^{\infty}(\R)$.
If $x$ is nondegenerate or of Morse-Bott type then this limit exists
for $R\to \infty$ and the asymptotic approach is
exponential.
\end{theorem}

Punctures of a finite energy holomorphic sphere can be either
positive or negative depending upon whether the image of the curve
lies in $(-\infty,0)\times M$ or $(0,\infty)\times M$ near the
puncture. In theorem $2$ the sign in the exponent is $+1$ in the
case of a convex end and $-1$ in the case of a concave end.

There is also a Fredholm theory for such curves. This was done in
\cite{hoff} in the nondegenerate case and in our degenerate
situation the result is again discussed in \cite{bor}, with the
proofs given in \cite{bort}.

Let $u$ be an embedded finite energy holomorphic sphere with
positive ends asymptotic to Reeb orbits
$\gamma^+_1$,...,$\gamma^+_{s^+}$ and negative ends asymptotic to
Reeb orbits $\gamma^-_1$,...,$\gamma^-_{s^-}$. We are interested
in the virtual dimension of the moduli space of finite energy
spheres containing $u$, modulo reparameterizations. This is the
index, $\mathrm{index}(u)$, of a certain Fredholm operator.

For generic choices of almost-complex structure $J$ satisfying
equation $(1)$ near the open ends, this index does indeed give the
dimension of the moduli space of finite energy holomorphic spheres
in a neighborhood of an embedded curve. As usual, virtual cycle
techniques must be employed to deal with multiply covered curves,
see the discussion in \cite{bort}. In this paper, all of the
finite energy curves we encounter will turn out to be embedded.
The theorem below can be found in \cite{bor} or \cite{bort}.

\begin{theorem}
With $u$ as above, the deformation index of $u$ is given by
$$\mathrm{index}(u)=-(2-s^+-s^-)+2c_1(TW)[u]+ \sum_{i=1}^{s^+} (\mu(\gamma^+_i)+\frac{1}{2}\mathrm{dim}(\gamma^+_i)) -\sum_{i=1}^{s^-} (\mu(\gamma^-_i)-\frac{1}{2}\mathrm{dim}(\gamma^-_i))$$
where $\mu(\gamma^{\pm}_i)$ is a generalized Conley-Zehnder index
defined in \cite{rs} and $\mathrm{dim}(\gamma^{\pm}_i)$ is the
dimension of the manifold of Reeb orbits containing
$\gamma^{\pm}_i$.
\end{theorem}

The moduli space containing a sphere $u$ will in general contain
spheres asymptotic to a different set of Reeb orbits. The
definitions of the Conley-Zehnder index and Chern class here are
given with respect to a fixed trivialization along the Reeb
orbits. More precisely, for each $i$ we choose a symplectic
trivialization of $\xi$ along $\gamma_i$. With respect to this
trivialization the Reeb flow gives a family of symplectic matrices
$d\eta_t \in \mathrm{Sp}(2, \Bbb R)$ for $0 \le t \le T$ where $T$
is the period of $\gamma_i$. We associate the index
$\mu(\gamma_i)$ to this family following \cite{rs}. Now, our
trivialization along the $\gamma_i$ naturally induces one of
$TW|_{\gamma_i}$ (thinking of $\gamma_i$ here as lying in
$(0,\infty) \times M$ or $(-\infty,0) \times M$) since
$T_{\gamma_i(t)}W=\xi_{\gamma_i(t)}\oplus \Bbb R \frac{\p}{\p t}
\oplus \Bbb R X$ and so of the complex line bundle $\bigwedge ^2
TW$ along the $\gamma_i$. This can be used to define
$c_1(TW)[u]=c_1(\bigwedge ^2 TW)[u]$ as follows. We choose a
section of $u^*\bigwedge ^2 TW$ which coincides with our
trivialization near the punctures and then count the numbers of
zeros with multiplicity.

Finally we state a compactness theorem which will be needed in the
course of our proof. Now let $(W, \omega, J)$ be a closed symplectic
manifold with a compatible almost-complex structure $J$. The
particular situation in which we will be interested is taking a limit
when we deform the almost-complex structure in the neighbourhood of
a contact-type hypersurface $\Sigma \subset W$, called stretching-the-neck.

A contact type hypersurface $\Sigma \subset W$ is an embedded contact
manifold with a contact form $\alpha$ such that $\omega|_{\Sigma}=d\alpha$.
This condition allows us to find a symplectic embedding of
$((-\epsilon,\epsilon)\times \Sigma,d(e^t\alpha))$ to a neighbourhood
$V$ of $\Sigma$ taking $\{0\}\times \Sigma$ onto $\Sigma$.
By perturbing $J$ near $V$ we may assume that it coincides
with the push-forward of an almost-complex structure given by formula
$(1)$ on $(-\epsilon,\epsilon)\times \Sigma$.

We suppose that $\Sigma$ divides $W$ into two symplectic manifolds $W_1$ and
$W_2$ such that $\Sigma$ is a convex boundary for $W_1$ and a concave
boundary for $W_2$. This means that the above embedding maps
$(-\epsilon,0)\times \Sigma$ into $W_1$ and $(0,\epsilon)\times \Sigma$
into $W_2$.

Following \cite{egh} and \cite{hwz}, we remove the tubular
neighbourhood $V$ of $\Sigma$ from $W$ and for each $N$ replace it
by gluing in a copy of $(-N,N) \times \Sigma$. We call the
resulting manifold $A_N$ and define an almost-complex structure
$J_N$ on $A_N$ by again using formula $(1)$ on $(-N,N) \times
\Sigma$ and letting $J_N=J$ elsewhere. This almost-complex
structure is compatible with a symplectic form $\omega_N$ on $A_N$
and $(A_N,\omega_N)$ is symplectomorphic to $(W,\omega)$ via a
symplectomorphism equal to the identity away from $V$.

We will need to study finite energy holomorphic spheres in three
associated noncompact symplectic manifolds with cylindrical ends.
Let $\tilde{W_1}$ be a completion of $W_1$ formed by gluing an end
symplectomorphic to $((0,\infty)\times \Sigma,d(e^t\alpha))$ and equipped with
a compatible almost-complex structure agreeing with $J$ on the
contact planes $\xi \subset T\Sigma$. Similarly define $\tilde{W_2}$
to be a completion of $W_2$ with end symplectomorphic to
$((-\infty,0)\times \Sigma,d(e^t\alpha))$ and equipped with a corresponding
compatible almost-complex structure. Third we have the
symplectization $S\Sigma =(\Bbb R\times \Sigma, d(e^t\alpha))$ which
also has a compatible almost-complex structure agreeing with $J$ on
the contact planes.

The following definitions and results are extracted from a more
detailed discussion in \cite{egh}, see also \cite{bor}. For a
proof see \cite{behwz} or \cite{bort}.

\begin{definition}
Let $(S,j)$ be a genus $0$ Riemann surface with nodes.
Then a level $k$ holomorphic map will consist of the following data:

(i) A labelling of the components of $S\setminus \{\mathrm{nodes}\}$ by
integers $\{1,...,k\}$ called levels such that two components sharing
a node have levels differing at most by $1$. Let $S_r$ be the union of
components of level $r$.

(ii) Finite energy holomorphic spheres $v_1:S_1\to \tilde{W_1}$,
$v_r:S_r\to S\Sigma$, $2\le r \le k-1$, and
$v_k:S_k\to \tilde{W_2}$. We require that each node shared by
$S_r$ and $S_{r+1}$ is a positive puncture for $v_r$ asymptotic to a
Reeb orbit $\gamma$ and a negative puncture for $v_{r+1}$ asymptotic to
the same Reeb orbit $\gamma$. Further $v_r$ should extend continuously
across each node within $S_r$.
\end{definition}

Suppose that $u_N:S^2\to (A_N,J_N)$ are a sequence of $J_N$
holomorphic curves where $S^2$ is the Riemann sphere with its
complex structure $i$. We suppose that the curves in the sequence
have bounded symplectic area. This is guaranteed if for instance
they lie in a fixed homology class, the situation we encounter in
this paper.

\begin{definition}
The sequence $u_N$ converges to a level $k$ holomorphic map from a
Riemann surface with nodes $(S,j)$ if there exist maps
$\phi_N:S^2\to S$ and sequences $t^r_N\in \Bbb R$, $r=2,...,k-1$, such that

(i) the $\phi_N$ are diffeomorphisms except that they may collapse circles
in $S^2$ to nodes in $S$, and $\phi_{N*}i \to j$ away from the nodes of $S$;

(ii) the sequences of maps $u_N\circ \phi_N^{-1}:S_1\to \tilde{W_1}$,
$u_N\circ \phi_N^{-1}+t^r_N:S_r\to S\Sigma$, $2\le r\le k-1$, and
$u_N\circ \phi_N^{-1}:S_k\to \tilde{W_2}$ converge in the $C^{\infty}$
topology to the corresponding maps $v_r$ on compact subsets of $S_r$.
\end{definition}

In the above definition, as is necessary we are identifying
$(-N,N)\times \Sigma \subset A_N$ with an increasing sequence of domains in $S\Sigma$,
$W_1\cup (-N,N)\times \Sigma \subset A_N$ with an increasing sequence of domains
in $\tilde{W_1}$ and $W_2\cup (-N,N)\times \Sigma \subset A_N$ with an increasing
sequence of domains in $\tilde{W_2}$.

\begin{theorem}
There exists a subsequence $N(i)$ of $N$ such that the sequence $u_{N(i)}$ converges to a level $k$
holomorphic map.
\end{theorem}

\section{Holomorphic Curves in $T^* S^2$}

Using the round metric on $S^2$, we can identify $T^* S^2$ with $TS^2$.
The aim of this section is to construct on $TS^2$ a convenient symplectic
form and almost-complex structure and describe possible foliations by
finite energy curves.

We will write $T^r S^2$ for the collection of tangent vectors of
length $r$. The pullback $\lambda$ of the Liouville form $pdq$
from $T^* S^2$ restricts to a contact form $\alpha$ on $T^1 S^2$
and the corresponding Reeb flow coincides with the geodesic flow.
In particular, it is periodic with period $2\pi$. We will call a
periodic orbit simple if its period is $2\pi$.

Denote by $\pi:TS^2\setminus S^2 \to T^1 S^2$ the projection
along the fibers from the complement of the zero-section to the
unit tangent bundle.
Let $\phi :TS^2\to [0,\infty)$ be a smooth increasing function
such that $\phi|_{T^r S^2}=r$ for $r\le 1$ and $\phi|_{T^r S^2}=e^r$ for $r$ large.
Then the $2$-form defined by $\omega=d(\phi \pi^* \alpha)$ extends
to a symplectic form on $TS^2$ equal to $d\lambda$ near the zero-section.
Globally, $(TS^2, \omega)$
is symplectomorphic to $(T^* S^2, d(pdq))$ via a symplectomorphism
fixing the zero-section.

We can think of $(TS^2,\omega)$ as one of our open symplectic manifolds
with a cylindrical end.
Observe that $SO(3)\equiv \mathrm{Isom}(S^2)$ acts
by differentials on $TS^2$. This action is by symplectomorphisms
preserving each $T^r S^2$.

We equip $(TS^2,\omega)$ with a compatible almost-complex
structure $J_0$ satisfying the following conditions. The
almost-complex structure $J_0$ should be invariant under the
action of $SO(3)$ on $TS^2$; the contact planes $\xi
=\mathrm{ker}(\alpha)$ on $T^r S^2$ should be invariant under
$J_0$; for $r$ sufficiently large $J_0$ should be invariant under
translation in the $r$ direction and $J_0(\frac{\p}{\p r})=X_r$.
For example, on the unit tangent bundle $J_0$ could be taken to be
the standard almost-complex structure mapping vertical tangent
vectors to their corresponding horizontal tangent vectors. For $r$
large, $J_0$ coincides with one of the standard almost-complex
structures on cylindrical ends given by formula $(1)$.

We can now study finite energy holomorphic spheres in $TS^2$. For
the moment we will assume that the almost-complex structure $J_0$
is suitably perturbed to an almost-complex structure $J$ such that
the linearization of the Fredholm operator from Theorem $3$ is
surjective and so its index gives the dimension of the
corresponding moduli space of finite energy spheres.

In the case of $TS^2$ there are no negative ends, so a finite
energy sphere $u$ has only positive asymptotic limits,
$\gamma_1$,...,$\gamma_s$. Suppose that $\gamma_i$ covers a simple
closed orbit $\mathrm{cov}(\gamma_i)$ times. We suppose that $u$
is embedded.

\begin{lemma} The dimension of the moduli space of finite energy
planes containing $u$ is given by
$$\mathrm{index}(u)=2(s-1)+\sum_{i=1}^{s} 2\mathrm{cov}(\gamma_i).$$
\end{lemma}

{\bf Proof} We choose a global trivialization of $\xi$ over $T^1
S^2$ of horizontal and vertical tangent vectors in $TS^2$. The
induced trivialization of $\bigwedge ^2 T(TS^2)$ over $T^1 S^2$
extends over all of $TS^2$ because of the existence of a global
splitting of $T(TS^2)$ into Lagrangian vertical and horizontal
subspaces. Thus in the formula of Theorem $3$ the Chern class term
will always be zero. Further, $\mathrm{dim}(\gamma_i)=2$ for all
$i$ since we have only one family of Reeb orbits. Thus to apply
Theorem $3$ it remains to compute the Conley-Zehnder index of an
orbit $\gamma$ with respect to this trivialization, say $\{H,V\}$
where $H$ and $V$ are unit horizontal and vertical vectors in
$\xi$ respectively. Given a vector $v\in \xi_{\gamma(0)}$, the
image under the Reeb flow $d \eta_t(v) = \frac{d}{ds}|_{s=0}
\gamma_s(t)$ where $\gamma_s$ is a family of closed orbits with
$\frac{d}{ds}|_{s=0} \gamma_s(0) =v$. Write $d\eta_t(v)=u_t H+w_tV
\in \xi_{\gamma(t)}$. Since all Reeb orbits $c$ in $T^1 S^2$
correspond to geodesics $\overline{c}$ in $S^2$ we observe that
$(u_t,w_t)=(J(t),J'(t))$ where $J(t)$ is the component of the
Jacobi field along $\overline{\gamma}$ corresponding to the
variation $\overline{\gamma_s}$. This Jacobi field is
perpendicular to $\overline{\gamma}$ since $d\eta_t(v) \in
\xi_{\gamma(t)}$. The Jacobi equation for the sphere is $J''+J=0$
and hence $(u_t,w_t)=(\cos(t)u_0 +\sin(t)w_0, -\sin(t)u_0 +
\cos(t)w_0)$ or equivalently $$d\eta_t =  \left( \begin{array}{cc}
\cos(t) & \sin(t) \\ -\sin(t) & \cos(t) \end{array} \right)  \in
\mathrm{Sp}(2,\Bbb R), 0 \le t \le 2\pi\mathrm{cov}(\gamma).$$ But
this path has Conley-Zehnder index $2\mathrm{cov}(\gamma)$ and the
result follows. \qed

Actually we can note that there are general formulas relating
Conley-Zehnder indices of closed Reeb orbits in unit tangent
bundles with Morse indices of the corresponding geodesics, see the
discussion in \cite{mohnke}.

From the lemma we see that $\mathrm{index}(u)\ge 2$ with equality
if and only if $u$ has a single puncture and simply covers a Reeb
orbit at that puncture.

Now, there are various possible foliations of $TS^2$ by finite
energy holomorphic spheres. For example, there exists a foliation
by finite energy planes in which all curves are asymptotic to a
cover of the same Reeb orbit (this was done in \cite{hind} in a
nondegenerate situation but the same foliation exists here). We
will obtain a foliation by a stretching-the-neck procedure. This
is described in section $4$ and we prove there that the foliation
is by finite energy planes asymptotic to simple Reeb orbits. Lemma
$10$ proves this and that there is a single plane asymptotic to
each such orbit. Meanwhile in this section we derive some
properties of such a foliation.

\begin{lemma}
Suppose that $TS^2$ is foliated by finite energy planes with a
single plane asymptotic to each simple Reeb orbit. Then the
intersection number of each plane with the zero-section is $\pm
1$.
\end{lemma}

{\bf Proof} Think of $TS^2$ as a neighborhood $U$ of
$\overline{\Delta}$ in $S^2 \times S^2$ with smooth boundary
$\Sigma=T^1 S^2$. Then $S^2 \times S^2 \setminus U$ is a disk
bundle over $\Delta$ and the boundaries of the disks are the Reeb
orbits in $\Sigma$. The asymptotic behaviour of our finite energy
planes allows us to compactify them to maps $(D,\p D)\to
(TS^2=U,\Sigma)$, and we can glue the boundaries to the disks in
the complement of $TS^2$ to obtain a foliation of $S^2 \times
S^2$, smooth at least away from $\Sigma$. As the spheres in the
foliation intersect $\Delta$ in a single point each and
necessarily have self-intersection number $0$, they lie in one of
the classes $[S^2 \times \mathrm{pt}]$ or $[\mathrm{pt} \times
S^2]$. In particular, the finite energy planes in $TS^2$ intersect
$\overline{\Delta}$ with intersection number $\pm 1$. \qed

Such a foliation of finite energy planes will be shown to exist
with respect to any regular almost-complex structure $J$ on
$TS^2$. We now consider a sequence of regular almost-complex
structures converging smoothly to the $SO(3)$ invariant structure
$J_0$. The author does not know whether or not this structure can
be assumed regular, that is, whether or not our index formula is
still valid.

In any case, following \cite{hwz}, by a similar procedure by which
we will later stretch-the-neck, using the compactness theorem we
can take a limit of a subsequence and find a $J_0$-holomorphic
finite energy plane through any point in $TS^2$. Taking a diagonal
sequence of almost-complex structures we can find disjoint and
embedded planes through a dense set of points. There is no
bubbling here since simple orbits already have minimal period and
$TS^2$ contains no closed holomorphic curves. The planes stay
disjoint and embedded in the limit by an application of positivity
of intersections, see \cite{md}. Since the limiting planes are
still asymptotic to simple Reeb orbits they are not multiply
covered. By a further limiting process we can include these planes
in a foliation of $TS^2$, and, again by positivity of
intersections, planes in this foliation exhaust the possible
limits of our finite energy planes in this subsequence up to
reparameterization.

\begin{lemma}
The $J_0$-holomorphic foliation of $TS^2$ is by finite energy
planes asymptotic to distinct simple orbits of the Reeb flow. Each
plane intersects the zero-section transversally in a single point.
\end{lemma}

{\bf Proof} Let $C$ be a plane in our foliation asymptotic to a
Reeb orbit $\gamma$ corresponding to a geodesic
$\overline{\gamma}$, and let $K$ be the $S^1$ subgroup of
$SO(3)=\mathrm{Isom}(S^2)$ which preserves $\overline{\gamma}$,
that is, the group of rotations about a perpendicular axis. For
the first part, we choose $L \subset SO(3)$ to be a small disk
through the identity which is transverse to $K$. Then by the
$SO(3)$ invariance of $J_0$, for any $l\in L$ the plane $l.C$ is
also a finite energy plane asymptotic to a Reeb orbit close to
$\gamma$. Similarly to the argument in Lemma $8$, we deduce that
all such finite energy planes are disjoint. To see this, we
observe that as the asymptotic limits are disjoint, so are the
disks we can glue in $S^2 \times S^2 \setminus U$ to obtain
spheres in $S^2 \times S^2$ as in Lemma $8$. These spheres have
intersection number $+1$ with $\Delta$ and nonnegative
self-intersection number by the positivity of intersections (the
spheres corresponding to $C$ and $l.C$ intersect only in $U$).
Together this implies a self-intersection number of $0$ and that
the planes are disjoint. Therefore the planes form a foliation of
a neighborhood of $C$. Now, any other finite energy plane $C'$ in
our foliation asymptotic to $\gamma$ but disjoint from $C$ must
intersect some of the curves $l.C$. But this contradicts
positivity of intersections since $C'$ must be homotopic to $C$
(fixing $\gamma$ in a suitable compactification).

For the second part, we also observe that $k.C$ must coincide with
$C$ for all $k \in K$ (for otherwise these planes would intersect
some $l.C$ giving a contradiction as above). The orbits of $K$ on
$S^2$ consist of a point $p$, the antipodal point $q$ to $p$, and
circles around $p$. Our plane cannot intersect $S^2$ in a circle
for this would imply the existence of a holomorphic disk in $TS^2$
with boundary on $S^2$, a contradiction to Stokes' Theorem since
$\omega$ is positive on holomorphic curves whereas its primitive
vanishes on the zero-section $S^2$. Therefore the plane can
intersect $S^2$ only at the points $p$ and $q$. Now, $K$ acts
transitively on $T_p S^2 \subset T_p(TS^2)$ and $T_q S^2 \subset
T_q(TS^2)$, so these intersections must be transversal (otherwise
the plane would be tangent to $S^2$ at $p$ or $q$, a contradiction
since $S^2$ is Lagrangian while embedded holomorphic curves are
symplectic). Thus each intersection is transversal and will
contribute $\pm 1$ to our intersection number. Since the total
intersection number is also $\pm 1$ we deduce that our finite
energy planes (which are embedded) intersect $S^2$ transversally
in a single point as claimed. \qed

We close this section by remarking that the transversal
intersection property will remain true for the regular
almost-complex structures in our subsequence which are
sufficiently close to $J_0$ by the smooth convergence ensured by
the compactness theorem.

\section{Stretching the neck}

Suppose that $L\subset S^2\times S^2$ is a Lagrangian submanifold
homologous to $\overline{\Delta}$. By Weinstein's theorem, a
sufficiently small neighbourhood $U$ of $L$ can be symplectically
embedded into $T^* S^2$, taking $L$ to the zero-section.
Let $T^{\le R} S^2$ denote the metric tube of radius $R$ inside
$TS^2$, $R$ large, and suppose that as in the previous section we have
constructed a symplectic form $\omega=d(\phi \pi^* \alpha)$ on $T^{\le R} S^2$
such that $\phi \pi^* \alpha$ is equal to the standard Liouville form near
the zero-section.

Now, for $\epsilon$ small enough, $(T^{\le R} S^2,\epsilon \omega)$ can
be symplectically embedded into $U$, again sending the zero-section
to $L$.

Let $\Sigma\subset U\subset S^2\times S^2$ be the boundary of a
tubular neighbourhood of $L\subset S^2\times S^2$ symplectomorphic
to $(T^{\le R} S^2,\epsilon \omega)$. We can push forward our
$SO(3)$-invariant almost-complex structure on $(T^{\le R}
S^2,\epsilon \omega)$ and extend it to a compatible almost-complex
structure on $S^2 \times S^2$. If necessary we will perturb this
almost-complex structure slightly such that the index of the
Cauchy-Riemann operator does indeed give the dimension of our
moduli spaces of holomorphic curves in all cases after we stretch
the neck.

As described in \cite{gr}, for any such compatible almost-complex
structure $J$ on $S^2\times S^2$ there exist two corresponding foliations
${\cal{F}}_1$ and ${\cal{F}}_2$ by $J$-holomorphic spheres in the
classes $[\mathrm{point}\times S^2]$ and $[S^2\times \mathrm{point}]$. Each sphere
in ${\cal{F}}_1$ intersects each sphere in ${\cal{F}}_2$ transversally in
a single point. We will now deform the almost-complex structure in
a neighbourhood of $\Sigma$,
stretching the neck as described in section $2$ to get a sequence $J_N$ of
almost-complex structures on $S^2\times S^2$.
We note that our almost-complex structure is already in the standard form
$(1)$ near $\Sigma$.
For each $N$ we have corresponding
$J_N$-holomorphic foliations ${\cal{F}}_1$ and ${\cal{F}}_2$.

We apply the compactness theorem from section $2$ as we take the
limit $N\to \infty$. Different reparameterizations of a suitable
subsequence converge to finite energy holomorphic curves from
punctured spheres into one of three symplectic manifolds, namely
$TS^2$ with our almost-complex structure, the symplectization
$S\Sigma$ of $\Sigma$ with its translation invariant
almost-complex structure and the complement of $U$ in $S^2\times
S^2$ completed with a cylindrical end symplectomorphic to
$((-\infty,0) \times \Sigma,d(e^t\alpha))$. The limit of a
sequence of spheres can be thought of as a tree in which the
vertices are finite energy spheres and the edges connect finite
energy spheres with the same Reeb orbits as asymptotic limits.
Choosing reparameterizations of spheres in ${\cal{F}}_1$ which
pass through a chosen point of $U$, after taking a subsequence of
$N\to \infty$, in the limit we can find a finite energy sphere
passing through any point of $TS^2$. Now, taking a diagonal
subsequence, we can find a collection of finite energy spheres in
the limit which pass through a dense set of points in $TS^2$. A
second limiting process as in \cite{hwz} can be used to find a
finite energy sphere through every point of $TS^2$. These
punctured spheres actually form a foliation and
reparameterizations of any converging sequence of
$J_N$-holomorphic spheres must converge to one of these finite
energy spheres. This follows from the nature of the convergence.
Any intersections or singular points would also be seen as
intersections amongst $J_N$-holomorphic curves in the
corresponding foliation ${\cal{F}}_1$, see for instance \cite{md}.
This does not however immediately exclude the possibility of the
image of a limiting curve arising as a branched cover. As is
required for section $3$, we want to show that this foliation is
by planes asymptotic to simple Reeb orbits.

As a remark, we observe that whether or not such a foliation
arises as a result of this limiting process, such a foliation is
necessarily present in $TS^2$ with our symmetric compatible
almost-complex structure. One way of seeing this would be to
replace $L$ with $\overline{\Delta}$ so that the whole arrangement
in $S^2 \times S^2$ is invariant under the action of $SO(3)$. Then
every Reeb orbit must be an asymptotic limit of a sphere in  our
foliation and thus the foliation is by planes invariant under
$1$-parameter subgroups and asymptotic to simple Reeb orbits.

Taking further diagonal subsequences, the limiting process also
gives a foliation by finite energy spheres in the (completed)
complement of $U$, say $W$. By studying this foliation and the
spheres in the symplectization $S\Sigma$ of $\Sigma$, we are able
to derive some infomation about our foliation of $TS^2$. Now,
fixing the trivialization along the Reeb orbits as before, the
Chern class term $c_1(T(S^2\times S^2))[u]$ in Theorem $3$ is now
equal to $2$ for the components of finite energy spheres $u$
mapping to $W$ in our homology class. This is because
$c_1(T(S^2\times S^2))$ gives $2$ when evaluated on spheres in the
foliation ${\cal{F}}_1$ but our trivialization on $TS^2$ (and
$S\Sigma$) is chosen such that it gives $0$ on the punctured
spheres in $TS^2$ and $S\Sigma$.

We can now observe that the component of the limiting holomorphic
map which has image in $W$ must be connected. For otherwise we
could find such a finite energy sphere $v$ in $W$ with Chern class
$c_1(T(S^2\times S^2))[v]<2$. This is a contradiction since such a
sphere could be glued to some planes in $TS^2$ to produce a
symplectic sphere in $S^2\times S^2$ of Chern class less than $2$.
This also implies immediately that such spheres in $W$ are not
multiple covers and therefore that our index formula is valid (as
the almost-complex structure was chosen generically).

The spheres $u$ in $W$ have only negative asymptotic ends, say
$\gamma_1,...,\gamma_s$, and so we find by Theorem $3$ and the
computation of the Conley-Zehnder indices in Lemma $7$ that the
dimension of the corresponding moduli space is given by
$$\mathrm{index}(u)=2(s+1)-\sum_{i=1}^{s} 2\mathrm{cov}(\gamma_i).$$
In particular, it is at most two, and equals two only if all of
the asymptotic limits simply cover the Reeb orbits. Therefore
generic simple Reeb orbits appear as negative asymptotic limits
for curves in the foliation of $W$. (Fixing a set of these limits
would necessarily give a moduli space of dimension less than two.)

In the symplectization $S\Sigma$ of $\Sigma$, the finite energy
spheres which appear in  our limits must have a single positive
puncture. (The maximum principle implies that there must be at
least one positive puncture.) This follows because the limit curve
has a single component in $W$, but as we are dealing only with
curves of genus $0$, different positive asymptotic limits of
curves in $S\Sigma$ could not be connected to the same component.
If this positive puncture is asymptotic to a simple Reeb orbit
then the curve must actually be a cylinder $\Bbb R \times \gamma$
over this orbit. This is because $\int u^* \pi^* d\alpha \ge 0$
for all curves in $S\Sigma$, where $\pi$ denotes the projection
onto $\Sigma$, and equality occurs if and only if the curve is a
cylinder. But by the asymptotic convergence to Reeb orbits, this
integral is just the difference between the periods of the
positive and negative asymptotic limits.

Similarly, the components of our limiting curve in $TS^2$ must
have a single positive asymptotic limit. Since for a generic curve
(that is, its component in $W$ passes through an open dense subset
in $W$) its component in $W$ has only simple negative asymptotic
limits and any components in $S\Sigma$ are cylinders, the
components of the curve in $TS^2$ are planes with a simple
asymptotic limit (and in particular are not multiple covers). We
find such finite energy planes asymptotic to an open dense set of
Reeb orbits and use this infomation to deduce our final lemma as
claimed in section $3$.

\begin{lemma}
The asymptotic limits of the curves in our foliation of $TS^2$ are
simple Reeb orbits, and there is a single curve asymptotic to each
Reeb orbit.
\end{lemma}

{\bf Proof} The moduli space of finite energy planes in $TS^2$
asymptotic to simple Reeb orbits does indeed have dimension $2$ by
the formula of Lemma $7$, and the part of the moduli space close
to a given plane consists of planes asymptotic to the nearby Reeb
orbits (since our almost-complex structure is assumed regular).
Again, such nearby planes are automatically disjoint and as in
Lemma $9$ we see that there is a single plane asymptotic to each
simple orbit. The moduli space of such planes is compact since
their asymptotic limits have minimal period and no further
bubbling is possible. Thus, including planes asymptotic to the
remaining Reeb orbits gives a foliation of $TS^2$ which is our
foliation as required. \qed

Hence we can now apply the results of section $3$ to deduce that
if our almost-complex structure is chosen sufficiently close to
$J_0$ each plane will intersect the zero-section $S^2$
transversally in a single point. Thus, after taking the
subsequence, for $N$ sufficiently large the $J_N$-holomorphic
spheres in the foliation ${\cal{F}}_1$ must also intersect $L$
transversally in a single point. Taking a further subsequence, we
can assume that the same is true for spheres in the foliation
${\cal{F}}_2$.

In conclusion, we have shown the existence of an almost-complex structure
$J_N$ on $S^2\times S^2$, tamed by the standard symplectic form, such that
that the transverse foliations ${\cal{F}}_1$ and ${\cal{F}}_2$ by holomorphic
spheres have the property that each leaf in each foliation intersects
$L$ transversally in a single point.

\section{Conclusion of proof}

We want to construct a Lagrangian isotopy between $L$ and $\overline{\Delta}$.
In fact, since the group of symplectomorphisms of $S^2\times S^2$ which
act trivially on homology is connected (by a result of Gromov \cite{gr} it
is homotopic to $SO(3)\times SO(3)$), it will suffice to construct
a symplectomorphism taking $L$ to $\overline{\Delta}$.

Corresponding to $L$ we have an almost-complex structure $J$ such that the
holomorphic curves in the corresponding foliations ${\cal{F}}_1$ and ${\cal{F}}_2$
intersect $L$ transversally in single points. Similarly, the standard
foliations $\mathrm{point}\times S^2$ and $S^2\times \mathrm{point}$ are holomorphic
for the standard split complex structure $J_0$ and each curve intersects
$\overline{\Delta}$ in a single point.

Now, there is a unique extension of any diffeomorphism $L\to
\overline{\Delta}$ to a diffeomorphism $\phi$ of $S^2 \times S^2$
sending the transverse foliations corresponding to $J$ onto those
corresponding to $J_0$. Provided that our initial diffeomorphism
is chosen to be orientation preserving, $\phi$ will preserve the
complex orientation on the leaves and act trivially on homology.
Both $\phi^{-1 *}\omega$ and $\omega$ itself are compatible with
the almost complex structure $\phi_* J$ (in the second case
because $\phi_* J$ preserves the foliations $\mathrm{point}\times
S^2$ and $S^2\times \mathrm{point}$ which are orthogonal with
respect to $\omega$) and so $\omega_t=(1-t)\phi^{-1 *}\omega
+t\omega$ is a family of cohomologous symplectic forms on $S^2
\times S^2$ with respect to which $\overline{\Delta}$ is
Lagrangian (the forms are clearly closed and they are
nondegenerate since each is compatible with $\phi_* J$).

Using Moser's method, we write $\omega_t=\phi^{-1 *}\omega +d\beta_t$.
Then $d\beta_t |_{\overline{\Delta}}=0$ which implies that
$\beta_t|_{\overline{\Delta}}=dh_t$ for some function
$h_t:\overline{\Delta}\to \Bbb R$.
We can extend the $h_t$ smoothly to functions on $S^2\times S^2$,
replace our $\beta_t$ by $\beta_t -dh_t$ and therefore assume that
$\beta_t|_{\overline{\Delta}}=0$ for all $t$.

Now let $X_t$ be the unique solution of $X_t \rfloor \omega_t=\frac{d\beta_t}{dt}$. Then ${\cal{L}}_{X_t} \omega_t =d(\frac{d\beta_t}{dt})=\frac{d\omega_t}{dt}$
and since $\frac{d\beta_t}{dt}|_{\overline{\Delta}}=0$ we deduce that
$X_t$ must be tangent to $\overline{\Delta}$ along $\overline{\Delta}$.
It follows that the time-$1$ flow $\psi$ of the time-dependent vector field
$X_t$ is a diffeomorphism preserving $\overline{\Delta}$ and such that
$\psi^* \omega = \phi^{-1 *}\omega$. Hence $\psi \circ \phi$ is our
required symplectomorphism.

\end{document}

%% file: newmacro.tex
\newcommand{\R}{{\Bbb R}}
\newcommand{\C}{{\Bbb C}}

\def\NABLA#1{{\mathop{\nabla\kern-.5ex\lower1ex\hbox{$#1$}}}}
\def\Nabla#1{\nabla\kern-.5ex{}_#1}

\newcommand{\p}{{\partial}}

\newcommand{\qed}{\hfill$\Box$\medskip}
\newtheorem{theorem}{Theorem} 

\newtheorem{lemma}[theorem]{Lemma}

\newtheorem{definition}[theorem]{Definition}

%% file: lag5.bbl
\begin{thebibliography} {99}
\bibitem{arnold} V. I. Arnold, First steps in symplectic topology,
{\it Russ. Math. Surv.}, 41 (6) (1986), 1-21.
\bibitem{bor} F. Bourgeois, A Morse-Bott approach to contact
homology, {\it Symplectic and contact topology: interactions and
perspectives (Toronto, ON/Montreal, QC, 2001)}, 55-77, Fields
Inst. Commun., 35, Amer. Math. Soc., Providence, RI, 2003.
\bibitem{bort} F. Bourgeois, A Morse-Bott approach to contact homology, dissertation, Stanford University, 2002.
\bibitem{behwz} F. Bourgeois, Y. Eliashberg, H. Hofer, K. Wysocki and E.
Zehnder, Compactness results in Symplectic Field Theory, preprint
AIM 2003-14, math.SG/0308183.
\bibitem{cz} C.C. Conley and E. Zehnder, Morse type index theory for flows and periodic solutions of Hamiltonian equations, {\it Comm. Pure Appl. Math.}, 37(1984), 283-299.
\bibitem{egh} Y. Eliashberg, A. Givental and H. Hofer, Introduction to symplectic field theory, GAFA 2000 (Tel Aviv, 1999), {\it Geom. Funct. Anal.}, 2000, Special Volume, Part II, 560--673.
\bibitem{ep} Y. Eliashberg and L. Polterovich, Unknottedness of Lagrangian surfaces in symplectic $4$-manifolds, {\it Internat. Math. Res. Notices}, 11(1993), 295-301.
\bibitem{survey} Y. Eliashberg and L. Polterovich, The problem of
Lagrangian knots in four-manifolds, {\it Geometric topology
(Athens, GA, 1993)}, 313-327, AMS/IP Stud. Adv. Math., 2.1, Amer.
Math. Soc., Providence, RI, 1997.
\bibitem{gr} M. Gromov, Pseudo-holomorphic curves in symplectic manifolds, {\it Inv. Math.}, 82(1985), 307-347.
\bibitem{hind} R. Hind, Holomorphic Filling of $\Bbb RP^3$, {\it Comm. in Contemp. Math.}, 2(2000), 349-363.
\bibitem{hof} H. Hofer, K. Wysocki and E. Zehnder, A characterisation of the tight three-sphere, {\it Duke Math. J.}, 81(1995), no.1, 159-226.
\bibitem{hofa} H. Hofer, K. Wysocki and E. Zehnder, Properties of pseudoholomorphic curves in symplectisations I: Asymptotics, {\it Ann. Inst. H. Poincar\'{e} Anal. Non Lineaire}, 13(1996), no.3, 337-379.
\bibitem{hofi} H. Hofer, K. Wysocki and E. Zehnder, Properties of pseudoholomorphic curves in symplectisations II: Embedding controls and algebraic invariants, {\it Geom. Funct. Anal.}, 5(1995), no.2, 337-379.
\bibitem{hoff} H. Hofer, K. Wysocki and E. Zehnder, Properties of pseudoholomorphic curves in symplectisations III: Fredholm theory, {\it Topics in nonlinear analysis}, 381-475, Prog. Nonlinear Differential Equations Appl., 35, Birkh\"{a}user, Basel, 1999.
\bibitem{hofd} H. Hofer, K. Wysocki and E. Zehnder, Properties of pseudoholomorphic curves in symplectizations IV: Asymptotics with degeneracies, {\it Contact and symplectic geometry (Cambridge, 1994)}, 78-117, Publ. Newton Inst., 8, Cambridge Univ. Press, Cambridge, 1996.
\bibitem{hwz} H. Hofer, K. Wysocki and E. Zehnder, Finite energy foliations of tight three-spheres and Hamiltonian dynamics, {\it Ann. of Math.}, 157(2003), 125-255.
\bibitem{md} D. McDuff, The local behaviour of $J$-holomorphic curves in almost complex $4$-manifolds, {\it J. Diff. Geom.}, 34(1991), 143-164.
\bibitem{mohnke} K. Mohnke, How to (symplectically) thread the eye
of a (Lagrangian) needle, preprint AIM 2002-1, math.SG/0106139.
\bibitem{rs} J. Robbin and D. Salamon, The Maslov index for paths, {\it Topology}, 32(1993), 827-844.
\bibitem{seidel} P. Seidel, Lagrangian two-spheres can be symplectically knotted, {\it J. Diff. Geom.}, 52(1999), 145-171.


\end{thebibliography}
